\documentclass[final]{siamart171218}

\usepackage{amsmath,amssymb,amsfonts,mathtools}
\usepackage{graphicx}
\usepackage{bm}

\newtheorem{prob}{Problem}[section]
\numberwithin{equation}{section}

\title{Optimal Convergence of Quadratic Not-a-Knot Spline Interpolation on Arbitrary Partitions}
\headers{Not-a-Knot Interpolation with Quadratic Splines}{A. Greene}

\author{Aussie Greene\thanks{Department of Civil and Environmental Engineering,\\ Vanderbilt University, Nashville, TN (aussie.m.greene@vanderbilt.edu).}}

\begin{document}
\maketitle

\begin{abstract}
The not-a-knot condition provides a simple spline interpolation scheme that requires no derivative data at the endpoints. Whereas the classical not-a-knot scheme using odd-degree splines selects knots that coincide with the data points, this approach does not extend to the even-degree case. 
In this paper we establish the well-posedness of the not-a-knot interpolation problem for even-degree splines, where knots are selected at  midpoints of data values.
We then show that the quadratic case, introduced in [M. J. D. Powell, {\it Approximation Theory and
Methods}, Cambridge University Press, 1981], achieves $\mathcal{O}(h^3)$ convergence for $C^3$ functions on arbitrary partitions without mesh-ratio restrictions. The proof derives a linear system defining the interpolating spline, normalizes the system to bound the matrix norm of its inverse using a Neumann-series argument, and establishes global error bounds.
\end{abstract}



\section{Introduction}

In this paper we consider the following interpolation problem:
    \begin{prob} \label{interp}
        Suppose we are given $a=t_1<\cdots<t_n=b$ and associated values $y_1,\dots,y_n$. Find a function $s$ defined on $[a,b]$ such that\begin{equation}
        s(t_i)=y_i, \quad i=1,\dots,n.\end{equation}
    \end{prob}
A common choice of function spaces for $s$ is $S_d^r(\Delta)$, denoting the space of univariate splines of degree $d$ and smoothness $r$ on a partition $\Delta$ of $[a,b]$. Splines provide desirable smoothness and convergence properties. The not-a-knot interpolation scheme using cubic splines was proposed by de Boor~\cite{deBoor78} for interpolation problems where no information is given on the derivative of the interpoland. This scheme was first shown to converge without mesh-ratio restrictions by de Boor~\cite{deBoor85}, with Beatson proposing an alternative proof shortly after~\cite{Beatson86}.

In this paper we consider an analogous even-degree not-a-knot interpolation scheme that selects an interior set of midpoints as knots. The quadratic case of this interpolation scheme has previously been proposed; see Eq.~(18.40) of \cite{Powell81}. In Sect.~21.2 of \cite{Powell81}, the author analyzes properties of two related quadratic interpolation schemes. The first scheme, where interpolation points are chosen at the midpoints of a predefined knot sequence, is proven to have a uniformly bounded operator norm in Theorem~21.2. In our case where the knots are chosen at midpoints of given data values, however, it is shown that the operator norm cannot be uniformly bounded for arbitrary distributions of interpolation points.

Nevertheless, in this paper, we show that this interpolation scheme has optimal order convergence with a constant independent of the mesh-ratio for interpolands $f\in C^3[a,b]$. Our proof follows a linear-algebraic argument  inspired by Beatson's approach~\cite{Beatson86} where we construct a system solving for the first derivatives of the spline at interpolation points and bound the matrix norm of its inverse to show convergence. In Section~\ref{sec:nak} we define the not-a-knot interpolation scheme for even-degree splines and prove its well-posedness. Section~\ref{sec:linsys} derives the linear system solving for the unique quadratic not-a-knot spline interpolant, and its convergence is proved in Section~\ref{sec:convergence}.




\section{Not-a-knot interpolation with even-degree splines} \label{sec:nak}

For finding a function $s$ in Problem \ref{interp} using not-a-knot interpolation schemes, we work with the space $S_d^{d-1}(\Delta)$. To construct $\Delta:=\{x_i\}_{i=1}^{k+2}$, the number of interior knots $k$ must be selected such that the dimension of the spline space is equal to the number of interpolation conditions $n$. This gives \begin{equation} \label{dim}
        n=\dim S_{d}^{d-1}(\Delta):=k+d+1.
    \end{equation} 

The classical not-a-knot interpolation scheme is well-defined for odd-degree splines. The knots $\Delta$ are simply chosen to correspond with the data points $t_1,\dots,t_n$, and then the first and last $(d-1)/2$ interior knots are removed; see Sect.~1.4.3 of \cite{Schumaker15}. This results in $k=n-d-1$ interior knots, which satisfies condition \eqref{dim}.

For even-degree splines, we instead define the set of interior knots as the midpoints of a selection of intervals produced by the data values $\{t_i\}_{i=1}^n$. We have \begin{equation}\begin{aligned} \label{knots1}
    x_1&=a,\\
    x_i&=\frac{t_{i+\mu-1}+t_{i+\mu}}{2},\quad i=2,\dots,k+1,\\
    x_{k+2}&=b,
\end{aligned}\end{equation}
where $k=n-d-1$ and $\mu=d/2$. In other words, setting the knots equal to the midpoints of $[t_i,t_{i+1}]$ for $i=1,\dots,n-1$ reduces $k$ by 1, and then removing the first and last $\mu$ midpoints reduces $k$ by $2\mu=d$. This gives $k=n-d-1$ for even-degree splines, satisfying condition \eqref{dim}.

\begin{proposition}
    Given an even degree $d=2\mu$ for $\mu\geq1$, $n\geq d+1$, a knot sequence $\Delta$ as defined in \eqref{knots1}, and arbitrary data values $y_1,\dots,y_n$, there exists a unique spline interpolant $s\in S_d^{d-1}(\Delta)$ that satisfies $s(t_i)=y_i$ for $i=1,\dots,n$.
\end{proposition}

\begin{proof}
    Define an extended knot sequence \begin{equation} \label{exknots}
        \boldsymbol{\tau}:=\{\underbrace{a,\dots,a}_{d+1},x_2,\dots,x_{k+1},\underbrace{b,\dots,b}_{d+1}\},
    \end{equation} and let $\{N_i^{d+1}\}_{i=1}^n$ denote the associated normalized B-splines of degree $d$. By Schoenberg-Whitney (see Theorem~4.61 of \cite{Schumaker07}), the observation matrix $M:=(N_j^{d+1}(t_i))_{i,j=1}^n$ is nonsingular if and only if its diagonal entries are positive. The repeated knots in \eqref{exknots} give $N_1^{d+1}(t_1)=N_{n}^{d+1}(t_n)=1$. Each B-spline $N_i^{d+1}$ has support $[\tau_i,\tau_{i+d+1}]$ for $i=2,\dots,n-1$. We have
    \begin{equation}
        \tau_i=\begin{cases}
            a, & 2\leq i\leq d+1,\\
            \frac{t_{i-\mu-1}+t_{i-\mu}}{2}, & d+2\leq i\leq n-1,
        \end{cases}
    \end{equation}
    \begin{equation}
        \tau_{i+d+1}=\begin{cases}
            \frac{t_{i+\mu-1}+t_{i+\mu}}{2}, & 2\leq i\leq n-d,\\
            b, & n-d+1\leq i\leq n-1.
        \end{cases}
    \end{equation}
    In all cases, we see that $\tau_i<t_i<\tau_{i+d+1}$. It follows that $N_i^{d+1}(t_i)>0$ for $i=2,\dots,n-1$. Since all diagonal entries of $M$ are positive, $M$ is nonsingular. Therefore there exists a unique spline interpolant $s\in S_d^{d-1}(\Delta)$ that satisfies the given conditions.
\end{proof}

\section{Derivation of linear system for $S^1_2(\Delta)$ interpolation} \label{sec:linsys}
    Given a spline $\tilde s\in S_d^r({\tilde \Delta})$ defined on a set of knots $\tilde\Delta:=\{x_i\}_{i=1}^{m}$, ``removing a knot'' $x_i$ from $\Delta$ is equivalent to keeping the knot and setting \begin{equation} \label {dd}
        s^{(d)}_-(x_i)=s_+^{(d)}(x_i),
    \end{equation} where $s_-$ and $s_+$ are the spline segments connected by $x_i$. Therefore given a spline $\tilde s\in S_d^r(\tilde\Delta)$ such that $\tilde s$ satisfies equation \eqref{dd} and a spline $s\in S_d^r(\Delta)$ for $\Delta:=\{\tilde\Delta\setminus x_i\}$, $\tilde s(t)=s(t)$ for all $t\in[a,b]$. We exploit this fact to produce a well-behaved linear system for finding an equivalent spline to the quadratic not-a-knot interpolant by using a modified set of knots.

    Working with the space $S_2^1(\Delta)$ results in a set of $k=n-3$ interior knots. These knots correspond to the midpoints of the data values barring the first and last intervals. Equivalently, we can solve for a spline $\tilde{s}\in S_2^1(\tilde{\Delta})$ where the interior knots also include the midpoint of the final interval $\frac{t_{n-1}+t_n}{2}$ with the condition that \begin{equation} \label{nak}
        s''_-\left(\frac{t_{n-1}+t_n}{2}\right)=s''_+\left(\frac{t_{n-1}+t_n}{2}\right).
    \end{equation}

\begin{figure} [H]
  \centering
  \includegraphics[width = .6\linewidth]{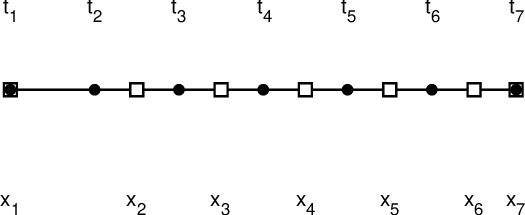}
  \caption{A set of data values $\{t_i\}_{i=1}^7$ and modified knots $\{x_i\}_{i=1}^7$.}
  \label{fig:knots}
\end{figure}

This produces a set $\tilde{\Delta}$ of $n$ knots with  \begin{equation}\begin{aligned}
    x_1&=a,\\
    x_i&=\frac{t_{i}+t_{i+1}}{2},\quad i=2,\dots,n-1,\\
    x_{n}&=b.
\end{aligned}\end{equation} Figure \ref{fig:knots} illustrates an example of $\tilde{\Delta}$ corresponding to a set of 7 uniformly spaced data points.

Let $\tilde s\in S_2^1(\tilde\Delta)$ be the spline that interpolates $t_i$ at $y_i$ for $i=1,\dots,n$. Define $h_i=t_{i+1}-t_i$ for $i=1,\dots,n-1$. Note that \begin{equation} \label{h2}
    t_{i+1}-x_i=\frac{h_i}{2},\quad i=2,\dots,n-1.
\end{equation} Let $S_i(t)$ be the segment of $\tilde s$ over the subinterval $[x_i,x_{i+1}]$ for $i=1,\dots,n-1$. Define $z_i=\tilde s'(t_i)$ for $i=1,\dots,n$. Because $S_1'(t)$ is linear and we have two interpolation conditions $t_1,t_2\in[x_1,x_2]$, as seen in Figure \ref{fig:knots}, we get \begin{equation} \label{S1d}
    S_1'(t)=\frac{z_1}{h_1}(t_2-t)+\frac{z_2}{h_1}(t-t_1).
\end{equation}
Now define $\sigma_i:=S_{i-1}'(x_i)$ for $i=2,\dots,n-1$. Using the $C^1$ continuity of $\tilde s$, the interpolation condition at $t_{i+1}\in[x_i,x_{i+1}]$, and \eqref{h2}, we obtain the following equations: \begin{equation} \label{Sid}
    S_i'(t)=\frac{2\sigma_i}{h_i}(t_{i+1}-t)+\frac{2z_{i+1}}{h_i}(t-x_i),\quad i=2,\dots,n-1.
\end{equation}

Following Beatson's approach~\cite{Beatson86}, we now seek to construct a linear system solving for $z=\begin{bmatrix}z_1&\cdots &z_n\end{bmatrix}^\top$ of the form \begin{equation} \label{Azb}
    Az=b,
\end{equation} where $A\in\mathbb{R}^{n\times n}$ and $b\in\mathbb{R}^n$. For our first equation, we integrate $S'_1(t)$ and use both interpolation conditions $t_1,t_2\in[x_1,x_2]$ to determine the constant of integration $C_1$. Integrating gives \begin{equation}
    S_1(t)=-\frac{z_1}{2h_1}(t_2-t)^2+\frac{z_2}{2h_1}(t-t_1)^2+C_1.
\end{equation} Then by plugging in the interpolation conditions, we get \begin{align}
     S_1(t_1)&=-\frac{z_1h_1}{2}+C_1=y_1,\\
      S_1(t_2)&=\frac{z_2h_1}{2}+C_1=y_2.
\end{align} Solving for $C_1$ gives our first equation \begin{equation} \label{A1}
z_1\frac{h_1}{2}+z_2\frac{h_1}{2}=y_2-y_1.\end{equation}

Integrating $S_i'(t)$ over the remaining intervals gives \begin{equation}
    S_i(t)=-\frac{\sigma_i}{h_i}(t_{i+1}-t)^2+\frac{z_{i+1}}{h_i}(t-x_i)^2+C_i,\quad i=2,\dots,n-1.
\end{equation} Using the interpolation conditions at $t_2,\dots,t_n$, we get \begin{equation}
    C_i=y_{i+1}-z_{i+1}\frac{h_i}{4},\quad i=2,\dots,n-1.
\end{equation} Now we derive the next $n-2$ equations using the $C^0$ continuity of $\tilde s$ at the interior knots $x_2,\dots,x_{n-1}$. Setting $S_1(x_2)=S_2(x_2)$ gives \begin{equation} \label{A2}
    z_1\left(-\frac{{h_2}^2}{8h_1}\right)+z_2\frac{4h_1h_2+{h_2}^2}{8h_1}+z_3\frac{h_2}{4}=y_3-y_2-\sigma_2\frac{h_2}{4}.
\end{equation} Then the remaining continuity conditions $S_i(x_{i+1})=S_{i+1}(x_{i+1})$ for $i=2,\dots,n-2$ are used to obtain the following $n-3$ equations: \begin{equation} \label{A3}
    z_{i+1}\frac{2h_ih_{i+1}+{h_{i+1}}^2}{4h_i}+z_{i+2}\frac{h_{i+1}}{4}=y_{i+2}-y_{i+1}+\sigma_i\frac{{h_{i+1}}^2}{4h_i}-\sigma_{i+1}\frac{{h_{i+1}}}{4},
\end{equation} for $i=2,\dots,n-2$.

To obtain our final equation, we enforce the not-a-knot condition \eqref{nak} upon the last interior knot. By deriving $S'_i(t)$ for $i=2,\dots,n-1$, we see that \begin{equation}
    S''_i(t)=-\frac{2\sigma_i}{h_i}+\frac{2z_{i+1}}{h_i},\quad i=2,\dots,n-2.
\end{equation} Setting $S''_{n-2}(x_{n-1})=S''_{n-1}(x_{n-1})$ yields the equation \begin{equation} \label{An}
    z_{n-1}\frac{1}{h_{n-2}}+z_{n}\left(-\frac{1}{h_{n-1}}\right)=\sigma_{n-2}\frac{1}{h_{n-2}}-\sigma_{n-1}\frac{1}{h_{n-1}}.
\end{equation}

Equations \eqref{A1}, \eqref{A2}, \eqref {A3}, and \eqref{An} produce a system of $n$ equations and $2n-2$ unknowns $z_1,\dots,z_n$ and $\sigma_2,\dots,\sigma_{n-1}$. Writing this system in the form \eqref{Azb} gives the banded matrix\begin{equation}
    A=\begin{bmatrix}
        \frac{h_1}{2}&\frac{h_1}{2}\\
        -\frac{{h_2}^2}{8h_1} & \frac{4h_1h_2+{h_2}^2}{8h_1} & \frac{h_2}{4}\\
        & 0 & \frac{2h_2h_3+{h_3}^2}{4h_2} & \frac{h_3}{4}\\
        & & \ddots & \ddots & \ddots \\
        & & & 0 & \frac{2h_{n-2}h_{n-1}+{h_{n-1}}^2}{4h_{n-2}}&\frac{h_{n-1}}{4}\\
        & & & & \frac{1}{h_{n-2}} & -\frac{1}{h_{n-1}}
    \end{bmatrix}
\end{equation}
    as well as the vector $b$ with components \begin{equation}\begin{aligned}
        b_1&=y_2-y_1,\\
        b_2&=y_3-y_2-\sigma_2\frac{h_2}{4},\\
        b_i&=y_{i+1}-y_i+\sigma_{i-1}\frac{{h_i}^2}{4h_{i-1}}-\sigma_i\frac{h_i}{4},\quad i=3,\dots,n-1,\\
        b_n&=\sigma_{n-2}\frac{1}{h_{n-2}}-\sigma_{n-1}\frac{1}{h_{n-1}}.
    \end{aligned}\end{equation}
    
    To eliminate $\sigma_2,\dots,\sigma_{n-1}$ from the system, we first perform the row operations
    \begin{enumerate}
        \item $R_n\leftarrow R_n-\left(\frac{4}{{h_{n-1}}^2}\right)R_{n-1}$,
        \item $R_i\leftarrow R_i+\left(\frac{{h_i}^2}{{h_{i-1}}^2}\right)R_{i-1},\quad i=n-1,n-2,\dots,3$.
    \end{enumerate}
    We then use \eqref{S1d} and \eqref{Sid} to derive \begin{equation}
        \sigma_2=-\frac{h_2}{2h_1}z_1+\frac{2h_1+h_2}{2h_1}z_2,
    \end{equation} \begin{equation}
        \sigma_i=-\frac{h_i}{h_{i-1}}\sigma_{i-1}+\frac{h_{i-1}+h_i}{h_{i-1}}z_i,\quad i=3,\dots,n-1,
    \end{equation} respectively. Substituting $\sigma_2,\dots,\sigma_{n-1}$ yields the system of equations \begin{equation} \label{Am}
    A=\begin{bmatrix}
        \frac{h_1}{2}&\frac{h_1}{2}\\
        -\frac{{h_2}^2}{4h_1} & \frac{3h_1h_2+{h_2}^2}{4h_1} & \frac{h_2}{4}\\
        & \frac{{h_3}^2}{4h_2} & \frac{3h_2h_3+3{h_3}^2}{4h_2} & \frac{h_3}{4}\\
        & & \ddots & \ddots & \ddots \\
        & & & \frac{{h_{n-1}}^2}{4h_{n-2}} & \frac{3h_{n-2}h_{n-1}+3{h_{n-1}}^2}{4h_{n-2}}&\frac{h_{n-1}}{4}\\
        & & & & -\frac{2}{h_{n-1}} & -\frac{2}{h_{n-1}}
    \end{bmatrix},
\end{equation} with \begin{equation} \label{bm}
        \begin{aligned}
            b_1&=y_2-y_1,\\
            b_2&=y_3-y_2,\\
            b_i&=y_{i+1}-y_i+\frac{{h_i}^2}{{h_{i-1}}^2}(y_i-y_{i-1}),\quad i=3,\dots,n-1,\\
            b_n&=-\frac{4}{{h_{n-1}}^2}(y_n-y_{n-1}).
        \end{aligned}
    \end{equation}

\section{Convergence for $C^3$ functions} \label{sec:convergence}
In this section we show that the not-a-knot interpolant with $s\in S_2^1(\Delta)$ has order three convergence for $C^3$ functions without mesh-ratio restrictions. We first perform another series of row operations in order to scale $A$ such that we can use the Neumann series technique to bound $\|A^{-1}\|_\infty$. Let $a_{ii}$ denote the $i$th diagonal element of $A$. The operations are as follows:
\begin{enumerate}
    \item $R_i\leftarrow R_i/a_{ii},\quad i=1,\dots,n$,
    \item $R_2\leftarrow \left(\frac{3h_1+h_2}{3h_1+3h_2}\right)R_2+\left(\frac{2h_2}{3h_1+3h_2}\right)R_1$,
    \item $R_1=\left(\frac{3h_1+4h_2}{3h_1+3h_2}\right)R_1-R_2$,
    \item $R_n\leftarrow \left(\frac{3h_{n-1}+4h_{n-2}}{3h_{n-1}+3h_{n-2}}\right)R_n-R_{n-1}$.
\end{enumerate}
We now obtain the matrix \begin{equation} \label{A}
    \hat A=\begin{bmatrix}
        1&\frac{h_2}{3h_1+3h_2} & -\frac{h_1}{3h_1+3h_2}\\
        \frac{h_2}{3h_1+3h_2} & 1 & \frac{h_1}{3h_1+3h_2}\\
         & \ddots & \ddots &\ddots\\
         & & \frac{h_{n-1}}{3h_{n-2}+3h_{n-1}} & 1 & \frac{h_{n-2}}{3h_{n-2}+3h_{n-1}}\\
         & & -\frac{h_{n-1}}{3h_{n-2}+3h_{n-1}} & \frac{h_{n-2}}{3h_{n-2}+3h_{n-1}} & 1
    \end{bmatrix},
\end{equation} with \begin{equation} \label{b}
        \begin{aligned}
            \hat b_1&=-\frac{4h_1(y_3-y_2)}{h_2(3h_1+3h_2)}+\frac{(6h_1+4h_2)(y_2-y_1)}{h_1(3h_1+3h_2)},\\
            \hat b_i&=\frac{4h_{i-1}(y_{i+1}-y_i)}{h_i(3h_{i-1}+3h_i)}+\frac{4h_i(y_{i}-y_{i-1})}{h_{i-1}(3h_{i-1}+3h_i)},\quad i=2,\dots,n-1,\\
            \hat b_n&=\frac{(6h_{n-1}+4h_{n-2})(y_n-y_{n-1})}{h_{n-1}(3h_{n-2}+3h_{n-1})}-\frac{4h_{n-1}(y_{n-1}-y_{n-2})}{h_{n-2}(3h_{n-2}+3h_{n-1})}.
        \end{aligned}
    \end{equation}

\begin{lemma} \label{bound}
    $\|\hat A^{-1}\|_\infty$ is bounded above by a constant that is independent from the mesh-ratio. 
\end{lemma}
\begin{proof}
    Let $B=\hat A-I$. Then for all rows $i=1,\dots,n$, \begin{align*}
    \sum_j^n |b_{ij}|&=\left|\pm\frac{h_{k+1}}{3h_k+3h_{k+1}}\right|+\left|\pm\frac{h_{k}}{3h_k+3h_{k+1}}\right|\\
    &=\frac{h_k+h_{k+1}}{3h_k+3h_{k+1}}\\
    &=\frac{1}{3}
\end{align*} for some $k$. It follows that $\|B\|_\infty=\frac{1}{3}$. Therefore by the Neumann series argument,\begin{equation}
    \|\hat A^{-1}\|_\infty=\|(I+B)^{-1}\|_\infty\leq\frac{1}{1-\|B\|_\infty}=\frac{3}{2}.
\end{equation}
\end{proof}

We will also need the following lemma.

\begin{lemma} \label{taylor}
    Given a function $f\in C^3$, we have that \begin{equation}
        \|\hat{A}f'(\textbf{t})-\hat{g}(f,\textbf{t})\|_{max} \leq K_1 \bar{h}^2\|f^{(3)}\|_\infty
    \end{equation} for some constant $K_1$ that is independent of the mesh-ratio.
\end{lemma}
\begin{proof}

    Define \begin{equation}r_i(x)=f(x)-q_i(x),\end{equation} where $q_i(x)$ denotes the quadratic Taylor polynomial of $f(x)$ about $t_i$. Since $s$ is a quadratic spline and the not-a-knot condition yields a unique interpolant to Problem \ref{interp}, we have that $s$ reproduces quadratic polynomials exactly. Hence, \begin{equation}
        \hat{A}q_i'(\textbf{t})=\hat{g}(q_i,\mathbf{t}).
    \end{equation} Then by linearity, we have that \begin{equation}
        (\hat{A}f'(\textbf{t})-\hat{g}(f,\mathbf{t}))_i=(\hat{A}r_i'(\textbf{t})-\hat{g}(r_i,\textbf{t}))_i.
    \end{equation} Taylor's theorem gives
    \begin{equation}
        |r_i(x)| \leq \frac{1}{6}|x-t_i|^{3}\|f^{(3)}\|_\infty,\quad |r'_i(x)| \leq \frac{1}{2}|x-t_i|^{2}\|f^{(3)}\|_\infty.
    \end{equation}

    We will also need
    \begin{align*}
        |r_1(t_3)-r_1(t_2)|&=\left| \int_{t_2}^{t_3} r'_1(x)\,dx\right|\\
        &\leq \int_{t_2}^{t_3} \left| r'_1(x)\right|\,dx\\
        &\leq \frac{1}{2}\|f^{(3)}\|_\infty \int_{t_2}^{t_3} |x-t_1|^2\,dx\\
        &\leq \frac{1}{2}\|f^{(3)}\|_\infty h_2(h_1+h_2)^2.
    \end{align*}


    Because $r_i(t_i)$, $r_i'(t_i)$ vanish and $h_i/(h_i+h_j)\leq 1$, we have
    \begin{align*}
        |(\hat{A}f'(\textbf{t})-\hat{b})_1|&\leq |r_1'(t_1)|+\frac{h_2}{3h_1+3h_2}|r_1'(t_2)|+\frac{h_1}{3h_1+3h_2}|r_1'(t_3)|\\
        &\quad+\frac{4h_1(|r_1(t_3)-r_1(t_2)|)}{h_2(3h_1+3h_2)}+\frac{(6h_1+4h_2)(|r_1(t_2)-r_1(t_1)|)}{h_1(3h_1+3h_2)}\\
        &\leq \left(\frac{1}{6}{h_1}^2 +\frac{1}{6}h_1(h_1+h_2)
        +\frac{2}{3}h_1(h_1+h_2)+\frac{1}{3}{h_1}^2\right)\|f^{(3)}\|_\infty\\
        &\leq\frac{13}{6}\bar{h}^2\|f^{(3)}\|_\infty.
    \end{align*}

    For rows $i=2,\dots,n-1$, we have
    \begin{align*}
        |(\hat{A}f'(\textbf{t})-\hat{b})_i| &\leq \frac{h_i}{3h_{i-1}+3h_i}|r_i'(t_{i-1})|+|r_i'(t_i)|+\frac{h_{i-1}}{3h_{i-1}+3h_{i}}|r'_i(t_{i+1})|\\
        &\quad+\frac{4h_{i-1}(|r_i(t_{i+1})-r_i(t_i)|)}{h_i(3h_{i-1}+3h_i)}+\frac{4h_i(|r_i(t_i)-r_i(t_{i-1})|)}{h_{i-1}(3h_{i-1}+3h_i)}\\
        &\leq \left(\frac{1}{6}{h_{i-1}}^2+\frac{1}{6}{h_i}^2+\frac{2}{9}{h_i}^2+\frac{2}{9}{h_{i-1}}^2\right) \|f^{(3)}\|_\infty\\
        &\leq \frac{7}{9}\bar{h}^2 \|f^{(3)}\|_\infty.
    \end{align*}

    Using a similar approach to the first row, the last row gives
    \begin{align*}
        |(\hat{A}f'(\textbf{t})-\hat{b})_n|&\leq \frac{h_{n-1}}{3h_{n-2}+3h_{n-1}}|r_n'(t_{n-2})|+\frac{h_{n-2}}{3h_{n-2}+3h_{n-1}}|r_n'(t_{n-1})|\\
        &\quad + |r'_n(t_n)| + \frac{(6h_{n-1}+4h_{n-2})(|r_n(t_n)-r_n(t_{n-1})|)}{h_{n-1}(3h_{n-2}+3h_{n-1})}\\
        &\quad+\frac{4h_{n-1}(|r_n(t_{n-1})-r_n(t_{n-2})|)}{h_{n-2}(3h_{n-2}+3h_{n-1})}\\
        &\leq \left( \frac{5}{6}h_{n-1}(h_{n-2}+h_{n-1})+\frac{1}{2}{h_{n-1}}^2\right)\|f^{(3)}\|_\infty\\
        &\leq \frac{13}{6}\bar{h}^2\|f^{(3)}\|_\infty.
    \end{align*}
    In view of these bounds, we choose $K_1$ to be the maximum of the three constants. This yields our desired result \begin{equation}
        \|\hat{A} f'(\textbf{t})-\hat{g}(f,\textbf{t})\|_{max}\leq\frac{13}{6}\bar{h}^2\|f^{(3)}\|_\infty.
    \end{equation}
\end{proof}

We now present a result on the error bounds of the interpolation scheme. We denote the $L^\infty$ function norm over $[a,b]$ and the $\ell^\infty$ vector norm by $\|\cdot\|_\infty$ and $\|\cdot\|_{\max}$ respectively.

\begin{theorem}
    Let $s\in S_2^1(\Delta)$ be the quadratic not-a-knot interpolant to Problem \ref{interp} with knots as defined in $\eqref{knots1}$. Define $\bar h:=\max\{t_{i+1}-t_i\}$ for $i=1,\dots,n-1$. Let $\textbf{t}\in\mathbb{R}^n$, $\textbf{x}\in\mathbb{R}^{n-3}$, and $\textbf{y}\in\mathbb{R}^n$ denote the vectors of data points $t_1,\dots,t_n$, interior knots $x_2,\dots,x_{n-2}$, and values $y_1,\dots,y_n$ respectively. Suppose $f\in C^3$ such that $f(t_i)=y_i$ for $i=1,\dots,n$. Then $s$ satisfies \begin{equation} \label{intcond}
            s(\textbf{t})=f(\textbf{t})\quad \text{and}\quad \hat As'(\textbf{t})= \hat g(\textbf{y},\textbf{t})
    \end{equation} for some $\hat A\in \mathbb{R}^{n\times n}$ such that $\|\hat A^{-1}\|_\infty$ is bounded and some $\hat g:\mathbb{R}^{n\times 2}\to \mathbb{R}^n$; and \begin{equation}
        \|f^{(k)}-s^{(k)}\|_\infty\leq K \bar h^{3-k}\|f^{(3)}\|_\infty
    \end{equation} for $k=0,1$ and some constant $K$ that is independent of the mesh-ratio.
\end{theorem}

\begin{proof}
    By construction, the system of equations $Az=b$ presented in \eqref{Am}, \eqref{bm} satisfy \eqref{intcond} where $b=g(\textbf{y},\textbf{t})$ for some $g\in \mathbb{R}^{n\times 2}\to \mathbb{R}^n$. Row operations may then be performed to obtain an equivalent system with matrix $\hat A$, as in \eqref{A}, for some vector $\hat{b}=\hat g(\textbf{y},\textbf{t})$. Lemma \ref{bound} shows that $\|\hat A^{-1}\|_\infty$ is bounded.

    Given a function $f\in C^3$, Lemma \ref{taylor} gives \begin{equation}
        \|\hat Af'(\textbf{t})-\hat g(f,\textbf{t})\|_{\max}\leq K_1 \bar h ^2 \|f^{(3)}\|_\infty,
    \end{equation} which implies \begin{equation} \|f'(\textbf{t})-s'(\textbf{t})\|_{\max}\leq K_1\|\hat A^{-1}\|_\infty \bar h^2 \|f^{(3)}\|_\infty.
    \end{equation}
    Let $\hat s'$ be the linear interpolant to $f'$ with knots $\{t_i\}_{i=1}^n$. By Theorem 6.15 in \cite{Schumaker07}, we have \begin{equation}
        \|f'-\hat s'\|_\infty \leq K_2 \bar h^2\|f^{(3)}\|_\infty,
    \end{equation} with $K_2=\frac{1}{8}$.
    
    Since $\hat{s}'-s'$ is piecewise linear with knots at $\textbf{t}$ and midpoints $\textbf{x}$,
    \begin{equation}\|\hat{s}'-s'\|_\infty=\max{\left\{\|\hat{s}'(\textbf{t})-s'(\textbf{t})\|_{\max},\|\hat{s}'(\textbf{x})-s'(\textbf{x})\|_{\max}\right\}}.\end{equation}
    Because $\hat{s}'$ interpolates $f'$, we have \begin{equation} \label{tmax}
        \|\hat{s}'(\textbf{t})-s'(\textbf{t})\|_{\max}=\|f'(\textbf{t})-s'(\textbf{t})\|_{\max}.
    \end{equation}
    Now consider a knot $x_i\in[t_i,t_{i+1}]$ for $i=2,\dots,n-2$. The trapezoidal rule gives
    \begin{equation} \label{trapezoidal}
        \int_{t_{i}}^{t_{i+1}}(\hat{s}'-s')(x)\,dx=\frac{h_i}{4}((\hat{s}'-s')(t_i)+2(\hat{s}'-s')(x_i)+(\hat{s}'-s')(t_{i+1})).
    \end{equation}
    Because $s$ interpolates $f$ at $t_i,t_{i+1}$, the fundamental theorem of calculus yields
    \begin{equation} \label{int}
                \int_{t_{i}}^{t_{i+1}}(\hat{s}'-s')(x)\,dx=\int_{t_{i}}^{t_{i+1}}(\hat{s}'-f')(x)\,dx.
    \end{equation}
    We combine \eqref{trapezoidal} and \eqref{int} then solve to get
    \begin{equation}
        (\hat{s}'-s')(x_i)=\frac{2}{h_i}\int_{t_i}^{t_{i+1}}(\hat{s}'-f')(x)\,dx-\frac{(\hat{s}'-s')(t_{i})+(\hat{s}'-s')(t_{i+1})}{2}.
    \end{equation}
    Hence,
    \begin{equation} \label{xmax}
        \|\hat{s}'(\textbf{x})-s'(\textbf{x})\|_{\max}\leq 2\|f'-\hat{s}'\|_\infty+\|f'(\textbf{t})-s'(\textbf{t})\|_{\max}.
    \end{equation}
    It follows from \eqref{tmax} and \eqref{xmax} that \begin{equation}
        \|\hat s'-s'\|_\infty\leq 2\|f'-\hat{s}'\|_\infty+\|f'(\textbf{t})-s'(\textbf{t})\|_{\max}.
    \end{equation}
    Then
    \begin{align*}
            \|f'-s'\|_\infty &\leq \|f'-\hat s'\|_\infty + \|\hat s'-s'\|_\infty\\
            &\leq 3\|f'-\hat s'\|_\infty + \|f'(\textbf{t})-s'(\textbf{t})\|_{\max}\\
            &\leq 3K_2 \bar h^2\|f^{(3)}\|_\infty + K_1\|\hat A^{-1}\|_\infty \bar h^2 \|f^{(3)}\|_\infty\\
            &=(3K_2+K_1\|\hat A^{-1}\|_\infty)\bar h^2\|f^{(3)}\|_\infty\\
            &\leq \frac{29}{8}\bar h^2\|f^{(3)}\|_\infty.
    \end{align*} 
    This proves the case for $k=1$. The case for $k=0$ follows from the fundamental theorem of calculus and the fact that $(f-s)(t_i)=0$ for $i=1,\dots,n-1$. For some $i$ and some $\hat t\in [t_i,t_{i+1}]$,
    \begin{align*}
        \|f-s\|_\infty&=|f(\hat t)-s(\hat t)|\\
        &=\left|\int_{t_i}^{\hat t} (f'-s')(\xi)\,d\xi\right|\\
        &\leq (t_{i+1}-t_{i})\|f'-s'\|_\infty\\
        &\leq \frac{29}{8}\bar h^3 \|f^{(3)}\|_\infty.
    \end{align*}
\end{proof}

We conclude that the quadratic not-a-knot interpolant achieves optimal $\mathcal{O}(\bar{h}^3)$ convergence over arbitrary partitions, without mesh-ratio restrictions.

\section{Acknowledgments}

The author would like to thank Dr. Larry L. Schumaker for suggesting the problem and funding his research during the author's appointment in the Department of Mathematics at Vanderbilt University.

\nocite{SharmaM66}

\bibliographystyle{siamplain}
\bibliography{references}

@book{deBoor78,
  author =	 "Carl de Boor",
  title =	 "A Practical Guide to Splines",
  publisher =	 "Springer--Verlag",
  year =	 "1978",
  address =	 "New York"
}

@Book{Schumaker15,
 author = {Larry L. Schumaker},
 title = {Spline Functions: Computational Methods},
 publisher = {SIAM},
 address = {Philadelphia},
 year = 2015
}

@article{Beatson86,
  author  = {R. K. Beatson},
  title   = {On the Convergence of Some Cubic Spline Interpolation Schemes},
  journal = {SIAM Journal on Numerical Analysis},
  volume  = {23},
  number  = {4},
  pages   = {903--912},
  year    = {1986},
}

@techreport{deBoor85,
  author      = {Carl de Boor},
  title       = {Convergence of Cubic Spline Interpolation with the Not-a-Knot Condition},
  institution = {Mathematics Research Center, University of Wisconsin--Madison},
  number      = {2876},
  year        = {1985},
}

@article{SharmaM66,
    author = {A. Sharma and A. Meir},
    title = {Degree of Approximation of Spline Interpolation},
    journal = {Journal of Mathematics and Mechanics},
    number = {5},
    volume = {15},
    year = {1966},
    pages = {759-767}
}

@Book{Schumaker07,
 author= {Larry L. Schumaker},
 title = {Spline Functions: Basic Theory},
 publisher = {SIAM},
 address = {Philadelphia},
 year = 2007
}

@book{Powell81,
    author    = {Powell, M. J. D.},
    title     = {Approximation Theory and Methods},
    publisher = {Cambridge University Press},
    address   = {Cambridge},
    year      = {1981},
    series    = {Cambridge Studies in Advanced Mathematics}
}

\end{document}